\documentclass[mathematics,article,submit,oneauthor,pdftex]{mdpi} 
\pdfoutput=1
\firstpage{1} 
\pubvolume{1}
\issuenum{1}
\articlenumber{0}
\pubyear{2021}
\copyrightyear{2020}
\datereceived{} 
\dateaccepted{} 
\datepublished{} 
\hreflink{https://doi.org/} % If needed use \linebreak
\Title{Mathematical modelling by help of category theory: models and relations between them}

\TitleCitation{Mathematical modelling by help of category theory: models and relations between them}

\Author{Dmitrii Legatiuk$^{1}$ \orcidA{}}

\AuthorNames{Dmitrii Legatiuk}

\AuthorCitation{Legatiuk, D.}
\address{%
$^{1}$ \quad Chair of applied mathematics, Bauhaus-Universit\"at Weimar, Germany; dmitrii.legatiuk@uni-weimar.de}

\corres{Correspondence: dmitrii.legatiuk@uni-weimar.de}

\abstract{The growing complexity of modern practical problems puts high demands on the mathematical modelling. Given that various models can be used for modelling one physical phenomenon, the role of model comparison and model choice becomes particularly important. Methods for model comparison and model choice typically used in practical applications nowadays are computation-based, and thus, time consuming and computationally costly. Therefore, it is necessary to develop other approaches for working abstractly, i.e. without computations, with mathematical models. The abstract description of mathematical models can be achieved by help of abstract mathematics, implying formalisation of models and relations between them. In this paper, a category theory-based approach to mathematical modelling is proposed. On this way, mathematical models are formalised in the language of categories, relations between the models are formally defined, as well as several practically relevant properties are introduced on the level of categories. Finally, an illustrative example is presented underlying how the category-theory based approach can be used in practice. Further, all constructions presented in this paper are also discussed from the modelling point of view by making explicit the link to concrete modelling scenarios.}

\keyword{category theory; mathematical modelling; abstraction; formal approaches; functors} 

\MSC{00A71; 06A75; 18B99; 18C10}
\begin{document}

\section{Introduction}

The rapid development of modern technologies naturally leads to higher demands for the mathematical modelling process, because practical problems nowadays require advanced coupled models. Moreover, typically several models can be used for modelling a given physical phenomenon, and thus, a model selection process must be made. Evidently, the model selection influences the quality of a final coupled model. In this regard, one of the most important tasks of a modeller is understanding the role of individual models in a complete coupled model, as well as studying how different models are related along with the practical meaning of this relation.\par
In engineering applications, various factors leading to reduction of the quality of the final coupled model are typically referred to as {\itshape uncertainties}. According to \cite{Oberkampf}, three types of uncertainties arising during the modelling process can be distinguished: (i) model inputs, (ii) numerical approximation, and (iii) model form. While the first two types can be identified and treated by help of computational and statistical methods, see for example \cite{Babuska,Keitel} and references therein, the third type requires an extra treatment. The model form uncertainty implies that a {\itshape conceptual modelling error} has been made, i.e. basic physical assumptions of models have been violated. Considering that the impact of such conceptual modelling errors on the whole modelling process is much more profound, it is necessary to develop tools towards addressing conceptual modelling errors.\par
Consideration of mathematical models based only on their physical assumptions, i.e. without considering a specific engineering example or performing computations with a model, requires tools of abstract mathematics. Several approaches to using abstract mathematics in applied mathematical modelling, such as graph theory \cite{Keitel}, abstract Hilbert spaces \cite{Dutailly_1,Dutailly_2}, abstract algebraic approach \cite{Nefzi,Legatiuk_2}, predicate logic \cite{Vassilyev_1,Vassilyev_2}, type theory \cite{Guerlebeck_1,Legatiuk_1}, and category theory \cite{Foley,Guerlebeck}, have been proposed in recent years. In this paper, we aim at revisiting and further developing the category theory-based modelling methodology introduced in \cite{Guerlebeck}. The motivation for using category theory for abstract description of mathematical models is based on several aspects: (i) the abstract nature of category theory allows description of very different objects and structures on common basis; (ii) a practical interpretation of abstract constructions provided by category theory-based modelling methodology is straightforward, and thus, the methodology can really be used in engineering practice; (iii) category theory naturally provides scaling possibilities implying that description of more sophisticated objects and structures can be done by using the same principles as descriptions of their individual parts; (iv) finally, various applications of category theory scattering from modelling of dynamical systems \cite{Behrisch} to ontological representation of knowledge \cite{Spivak} presented in recent years indicate that advantages of category theory are seen and accepted now not only by mathematicians, but also by people interested in applications.\par
As we have already mentioned, the category theory-based modelling methodology discussed in this paper has been originally proposed in \cite{Guerlebeck}. After publishing this work, several new ideas on categorical modelling methodology providing a deeper understanding of mathematical models and modelling process have appeared in recent years. Therefore, it is necessary to revise ideas presented in \cite{Guerlebeck} with new results and more refined categorical constructions. Moreover, it is worth to mention, that the use of category theory-based modelling methodology for analysis of models appearing in real-world engineering problems from the field of aeroelastic analysis of bridges has been presented in \cite{Kavrakov}. This work indicated practical advantages of using category theory for modelling purposes. To this end, the category theory-based modelling methodology presented in this paper aims at a consistent description of mathematical models and relations between them in the language of category theory. For the sake of clarity, we focus in this paper only on individual mathematical models, while coupled models will be treated in future research using results from the current paper as a basis.\par 
Abstract categorical descriptions of mathematical models requires at first defining {\itshape universal properties of models}, which are properties shared by models in general, i.e. independent on a particular problem of an engineering field. If a universal model property is defined, then all categorical constructions used in one specific modelling application can be directly transferred to another field. Thus, we will start our construction with defining such a universal model property which is common for all models. Moreover, the main goal is to keep track of real physical and engineering interpretations of the constructions introduced in the category theory-based modelling methodology. The paper is organised as follows: Section~\ref{Section:categories_models} presents a general structure of categories of mathematical models together with a detailed discussion on practical interpretation of the introduced definition; after that, relations between mathematical models are discussed in Section~\ref{Section:relations}; Section~\ref{Section:convertible} formalises the problem of having different formulations of the same mathematical model by introducing the notion of {\itshape convertible mathematical models}; Section~\ref{Section:example} provides an illustrative example how categorical constructions introduced in the previous sections can be used for comparison and analysis of models. Finally, in Section~\ref{Section:conclusions} we discuss a universal arrow in the framework of category theory-based modelling methodology, as well as establish a connection to an abstract algebraic approach, after we draw conclusions and discuss shortly the scope of future work. For making the paper self-contained, some basic definitions from category theory are presented in the Appendix.\par

\section{Categories of mathematical models}\label{Section:categories_models}

Before starting with categorical constructions, it is important to underline, that models used in practice can be generally classified in two types:
\begin{itemize}
\item {\itshape physics-based models} -- models which are based on mathematical formalisations of physical laws and assumptions;
\item {\itshape data-driven models} -- models which are based on representations of data, e.g. results of experiments or measurements obtained from a monitoring system.
\end{itemize}
This paper deals with physics-based models, which are referred to simply as {\itshape mathematical models}, because this type of models is typically implied by the term {\itshape mathematical modelling}. Moreover, because mathematical models are based on physical assumptions formalised by help of mathematical expressions, they provide a richer basis for abstract considerations, compared to data-driven models, which are very often black-box models not relying on any physical assumptions.\par
We start our construction with the introduction of {\itshape concrete} categories $\mathbf{Model}_{i}$, $i=1,2,\ldots$, which are associated with mathematical models used to describe a certain physical phenomenon, such as, for example, models of elasticity theory or heat conduction. The term \textquotedblleft associated\textquotedblright\, has been used, because, strictly speaking, the objects of categories $\mathbf{Model}_{i}$, $i=1,2,\ldots$ are not mathematical models themself, but rather sets of basic physical assumptions on which the corresponding mathematical models are created. However, to keep notations short and transparent, we will refer to these categories simply as to {\itshape categories of mathematical models}. The following definition introduces basic structure of these categories:
\begin{Definition}[Category of mathematical models]\label{Definition:category_of_mathematical_models}
Let $\mathbf{Model}_{1}$ be a category of mathematical models describing a given physical phenomenon. Then for all objects of $\mathbf{Model}_{1}$ the following assumptions hold:
\begin{itemize}
\item[(i)] each object is a finite non-empty set -- {\itshape set of assumptions} of a mathematical model, denoted by $\mathbf{Set}_{\mathrm{A}}$, where $\mathrm{A}$ is the corresponding mathematical model;
\item[(ii)] morphisms (arrows) are relations between these sets;
\item[(iii)] for each set of assumptions and its corresponding model exists a mapping

\begin{equation*}
\mathbf{Set}_{\mathrm{A}} \stackrel{S}{\mapsto} \mathrm{A};
\end{equation*}

\item[(iv)] all objects are related to mathematical models acting in the same physical dimension.
\end{itemize}
\end{Definition}
Let us now provide some motivation from the modelling perspective and comments for the assumptions used in this definition:
\begin{itemize}
\item {\bfseries Assumption (i)}. This assumption comes naturally from the modelling process: a mathematical model is created to describe a certain physical phenomenon or process, and evidently, it is possible only if physical background of the phenomenon or process is clearly stated, i.e. assumptions to be satisfied by the model are formulated. Moreover, for a stronger distinction between different mathematical models, the set of assumptions is understood in a broader sense: not only basic physical assumptions are listed, but all further modifications and simplifications of the model, such as for example a linearisation of original equations, are also elements of the set of assumptions. The requirements for the set of assumptions to be finite comes from the fact that no model possess an infinite set of physical assumptions. Therefore, consideration of more general sets is not necessary.\\[5pt] 
It is also important to remark that having finite sets as objects in the category is one possible way to approach mathematical models. Alternatively, one could think of working directly with mathematical expressions (equations) representing the models. However, in this case it will be more difficult to distinguish models, since the same set of assumptions can be formalised differently in terms of final equations, as we will see in Section~\ref{Section:convertible}.
\item {\bfseries Assumption (ii)}. This assumption, in fact, introduces the structure of categories of mathematical models. The main point here is that instead of working with discrete categories, it is beneficial to study more elaborated structure. Since the objects in categories of mathematical models are sets, it is natural to use relations between sets as morphism in the categories. We will make these relations more specific in Section~\ref{Section:relations}.
\item {\bfseries Assumption (iii)}. This assumption formally describes the process of obtaining the final form of a model, e.g. differential or integral equation, from basic physical assumptions. In this case, mapping $S$ is, in fact, a {\itshape formalisation process} consisting in writing basic physical assumptions in terms of mathematical expressions, which constitute a mathematical model in the end of the formalisation process. Naturally, the formalisation process can be done by different means and approaches, for example first ideas on using type theory to describe the formalisation process towards detecting conceptual modelling errors have been presented in \cite{Guerlebeck_1,Legatiuk_1}.\\[5pt]
We also would like to remark, that originally, mapping $S$ has been called {\itshape invertible} in \cite{Guerlebeck}. The invertibility in this case means, that set of assumptions can be uniquely reconstructed from the final form of a model. Although that such a reconstruction is theoretically indeed possible, it is generally not unique. Even if we consider the following canonical parabolic equation

\begin{equation*}
u_{t} = a^{2}u_{xx},
\end{equation*}

then without extra context it cannot be decided if this is a heat equation or a diffusion equation. Therefore, the invertibility of a mapping $S$ has been dropped from  Definition~\ref{Definition:category_of_mathematical_models}.
\item {\bfseries Assumption (iv)}. This assumption ensures that we do not treat equally models from different dimensions.
\end{itemize}
It is also important to mention that according to Definition~\ref{Definition:category_of_mathematical_models}, models with different parameters, e.g. material constants, will be corresponded to the same set of assumptions. For example, if we consider the set of assumptions leading to the Lam\'e equation (partial differential equation with constant coefficients), then it is clear that infinite number of constant coefficients exists, but all these specific models are originated from the same set of assumptions. In general, models originating from the same set of assumptions, but having different material parameters are just particular instance of a general set of assumptions. This fact is particularly important for engineering applications, where stochasticity of material parameters in deterministic models is often considered as stochastic modelling. However, as we discussed above, the stochasticity only in material parameters does not change basic modelling assumptions, because the fact that a constant is chosen according to a certain probability law does not principally affect the assumption of having constant coefficients. In contrast, modelling of physical process by help of stochastic partial differential equations is based on completely different modelling assumptions, see for example \cite{Holden}, and therefore, should not be put together with \textquotedblleft classical\textquotedblright\, mathematical models.\par

\section{Relations between mathematical models}\label{Section:relations}

This section is devoted to defining relations between sets of assumptions, which are objects in categories of mathematical models, as introduced in Definition~\ref{Definition:category_of_mathematical_models}. The main requirement for such relations is that their must define a {\itshape universal model property}, which is independent on a specific problem, meaning that boundary or initial conditions (but not coupling/transmission conditions!) should not have influence on the model property. For satisfying this requirement, the comparison of mathematical models by help of universal model property called {\itshape model complexity} is proposed \cite{Guerlebeck}:
\begin{Definition}[Complexity of mathematical models]\label{Defition:complexity}
Let $\mathrm{A}$ and $\mathrm{B}$ be mathematical models in a category $\mathbf{Model}_{1}$. We say that model $\mathrm{A}$ has higher complexity than model $\mathrm{B}$ if and only if $\mathbf{Set}_{\mathrm{A}}\subset\mathbf{Set}_{\mathrm{B}}$, but $\mathbf{Set}_{\mathrm{B}}\not\subset \mathbf{Set}_{\mathrm{A}}$. Consequently, two models are called equal, in the sense of complexity, iff $\mathbf{Set}_{\mathrm{B}}=\mathbf{Set}_{\mathrm{A}}$.
\end{Definition}
The model complexity in this definition is defined relatively, since we do not describe it explicitly. From the point of view of physics, model complexity reflects the fact that a model which has less assumptions provides a more accurate description of a physical phenomenon under consideration. Thus, the model complexity is a {\itshape relative quality measure} of how good a mathematical model represents a given physical phenomenon. The relativity in the measure comes from the fact, that any comparison needs at least two objects, and one model cannot be assessed with respect to its ability represent the corresponding physical process, otherwise that would imply that the exact representation of the physical process is known {\itshape a priori}.\par
It is important to underline, that the notion of model complexity proposed in Definition~\ref{Defition:complexity} is neither related to the notion of complexity of an algorithm, nor to the notion of complexity used for statistical models, where the number of parameters is typically served as complexity measure. The advantage of the notion of model complexity introduced in Definition~\ref{Defition:complexity} is the fact that it does not depend on specific boundary or initial conditions, since typically basic model assumptions are not influenced by them. Nonetheless, if boundary conditions are essential for basic model assumptions, e.g. singular boundary conditions, then they will be automatically listed in the corresponding set of assumptions, since such boundary conditions are critical for describing the physical process. Thus, the model complexity introduced in Definition~\ref{Defition:complexity} is a universal model property.\par
Additionally, Definition~\ref{Defition:complexity} might sound a bit counterintuitive, since it states that a model satisfying less modelling assumption is more complex, and not of higher {\itshape simplicity}, as it could be expected as well. In fact, both points of view on the complexity are possible, and differ only in the general understanding of modelling assumptions. Definition~\ref{Defition:complexity} is based on the idea that modelling assumptions act as {\itshape restrictions} for a model, and thus, implying that a model with less modelling assumptions is more general. Nonetheless, another perspective on the notion of model complexity still can be considered, which would reflect the opposite point of view that model assumptions are not restrictions, but rather {\itshape generalisations} of models. This discussion is also directly related to the following important remark:
\begin{Remark}
Sets of assumptions introduced in  Definition~\ref{Section:categories_models} are assumed to be written by help of a natural language. Although intuitively it is clear how to formulate these sets, as well as how to compare them in the sense of model complexity, from the formal perspective it is not so straightforward. In fact, a formal comparison of sets of assumptions written in a natural language can be done only by help of a detailed semantic analysis of these sentences, and only after that, sentences, and hence sets of assumptions, can be rigorously compared. As a possible way around this problem, stricter rules on formulating sets of assumptions might be imposed. In that case, a kind of basic \textquotedblleft alphabet\textquotedblright\, containing allowed expressions and symbols could be introduced. Moreover, perhaps a combination of a natural language and mathematical expressions complemented by strict rules could be a suitable option. Different possibilities to address the problem of a rigorous comparison of sets of assumptions will be studied in future work.
\end{Remark}\par
From the point of relational algebra, model complexity is a binary relation in a category of mathematical models. Hence, the objects in categories of mathematical models can be ordered by using model complexity. However, the ordering of objects defined by model complexity is only partial, and not total, since examples of mathematical models which should belong to the same category but cannot be ordered according to Definition \ref{Defition:complexity} can be easily found, see for example aerodynamic models used in bridge engineering \cite{Kavrakov}. Naturally, in some cases mathematical models can constitute a category with totally ordered objects. To have a clear distinction between categories with partial and total ordering of objects, we introduce the following definition \cite{Kavrakov}:
\begin{Definition}\label{Definition:partially_totally_categories}
Let $\mathbf{Model}_{1}$ be a category of mathematical models in which $n$ objects $\mathbf{Set}_{A_{j}}$, $j=1,\ldots,n$ can be ordered according to Definition \ref{Defition:complexity} as follows
\begin{equation*}
\mathbf{Set}_{A_{i}}\subset \mathbf{Set}_{A_{j}}, \mbox{ for } i<j\leq n.
\end{equation*}

Moreover, let $X$ be the set of all modelling assumptions used in this category. Then category $\mathbf{Model}_{1}$ contains totally ordered objects, and therefore is associated with totally ordered models, iff
\begin{equation*}
X = \mathbf{Set}_{\mathrm{A}_{1}}\cup\mathbf{Set}_{\mathrm{A}_{2}}\cup\ldots\cup\mathbf{Set}_{\mathrm{A}_{n}}, \mbox{ and } \mathbf{Set}_{\mathrm{A}_{n}}=X,
\end{equation*}
otherwise, the category $\mathbf{Model}_{1}$ contains partially ordered objects corresponding to partially ordered models.
\end{Definition}
As a direct consequence of this definition we have the following corollary:
\begin{Corollary}\label{Corollary:simplest_most_complex}
In a totally ordered category $\mathbf{Model}_{1}$ with $n$ objects always exist two unique objects:
\begin{itemize}
\item object $\mathbf{Set}_{\mathrm{A}_{1}}$ satisfying $\mathbf{Set}_{\mathrm{A}_{1}} \subset \mathbf{Set}_{\mathrm{A}_{i}}$ $\forall i=2,\ldots,n$, which is called the most complex object, and the associated model $\mathrm{A}_{1}$ is called the most complex model;
\item object $\mathbf{Set}_{\mathrm{A}_{n}}$ satisfying $\mathbf{Set}_{\mathrm{A}_{n}}=\mathbf{Set}_{\mathrm{A}_{1}}\cup\mathbf{Set}_{\mathrm{A}_{2}}\cup\ldots\cup\mathbf{Set}_{\mathrm{A}_{n}}$, which is called the the simplest object element, and the associated model $\mathrm{A}_{n}$ is called the simplest model.
\end{itemize}
\end{Corollary}\par
It is worth to mention, that in the framework of introduced modelling formalism, the most complex object and the simplest object are, in fact, {\itshape initial object} and {\itshape terminal object} in categories of mathematical models, respectively. Note that, although categories of mathematical models have finite sets as objects, the initial and terminal objects are different to the ones in the classical category $\mathbf{Sets}$, where these are given by the empty set and one-element set, correspondingly. The difference comes precisely from the modelling background of our categories, since while formally it is still possible to consider the empty and one-element sets as sets of assumptions of some (fictitious) models, it does not make sense from the modelling perspective.\par
The proof of Corollary~\ref{Corollary:simplest_most_complex} is straightforward, and we only would like to mention, that uniqueness of objects $\mathbf{Set}_{\mathrm{A}_{1}}$ and $\mathbf{Set}_{\mathrm{A}_{n}}$ follows immediately from Definition~\ref{Defition:complexity} and from the fact that a totally ordered category is considered. The situation is trickier in the case of partially ordered categories:
\begin{Proposition}\label{Proposition:partially_simplest_complex}
For a partially ordered category $\mathbf{Model}_{1}$ with $n$ objects one of the following statements holds:
\begin{itemize}
\item[(i)] the most complex object $\mathbf{Set}_{\mathrm{A}_{1}}$ and the simplest object $\mathbf{Set}_{\mathrm{A}_{n}}$ do not exist;
\item[(ii)] the most complex object $\mathbf{Set}_{\mathrm{A}_{1}}$ exists, while the simplest object $\mathbf{Set}_{\mathrm{A}_{n}}$ does not exist;
\item[(iii)] the most complex object $\mathbf{Set}_{\mathrm{A}_{1}}$ does not exist, while the simplest object $\mathbf{Set}_{\mathrm{A}_{n}}$ exists;
\item[(iv)] the most complex object $\mathbf{Set}_{\mathrm{A}_{1}}$ and the simplest object $\mathbf{Set}_{\mathrm{A}_{n}}$ exist simultaneously.
\end{itemize}
\end{Proposition}
\begin{proof}
We prove this proposition by straightforwardly constructing corresponding structures of partially ordered categories. We start the proof by proving cases (ii) and (iii) at first, since the proof of case (i) will be based on cases (ii) and (iii), and finally we will prove case (iv). We consider a category with one object $\mathbf{Set}_{\mathrm{A}_{1}}$, and the rest objects we construct explicitly from $\mathbf{Set}_{\mathrm{A}_{1}}$. Without loss of generality we assume that $\mathbf{Set}_{\mathrm{A}_{1}}$ contains at least one element, which will be denoted by $\mathrm{A}_{1}^{(1)}$. The objects $\mathbf{Set}_{\mathrm{A}_{2}}$ and $\mathbf{Set}_{\mathrm{A}_{3}}$ are then constructed from $\mathbf{Set}_{\mathrm{A}_{1}}$ by adding different elements $\mathrm{A}_{1}^{(2)}$ and $\mathrm{A}_{1}^{(3)}$ to $\mathbf{Set}_{\mathrm{A}_{1}}$, correspondingly, i.e. we obtain new sets of assumptions by adding two different assumptions. This construction is shown by the diagram

\begin{equation*}
\begin{tikzpicture}[scale=1,cap=round]
\tikzstyle{information text}=[rounded corners,fill=red!10,inner sep=1ex]
\begin{scope}[>=latex]
    \path (0,0) node (a) {$\mathbf{Set}_{\mathrm{A}_{1}} = \left\{\mathrm{A}_{1}^{(1)} \right\}$};
    \path (-2,-1.5) node (b) {$\left\{\mathrm{A}_{1}^{(1)}, \mathrm{A}_{1}^{(2)} \right\} =\mathbf{Set}_{\mathrm{A}_{2}}$};
    \path (2,-1.5) node (c) {$\mathbf{Set}_{\mathrm{A}_{3}} = \left\{\mathrm{A}_{1}^{(1)}, \mathrm{A}_{1}^{(3)} \right\}$};
    \draw[->,line width=1.0pt] (a) -- (b);
    \draw[->,line width=1.0pt] (a) -- (c); 
  \end{scope}
\end{tikzpicture}
\end{equation*}
implying that $\mathbf{Set}_{\mathrm{A}_{1}}\subset\mathbf{Set}_{\mathrm{A}_{2}}$ and $\mathbf{Set}_{\mathrm{A}_{1}}\subset\mathbf{Set}_{\mathrm{A}_{3}}$, but $\mathbf{Set}_{\mathrm{A}_{2}}$ and $\mathbf{Set}_{\mathrm{A}_{2}}$ are not related. Thus, $\mathbf{Set}_{\mathrm{A}_{1}}$ is the most complex object in this category, but no the simplest object exists. Thus, the case (ii) is proved.\par 
The proof of case (iii) is analogues to case (ii), where only instead of adding extra assumptions, we remove different assumptions from the initial set. Thus, for simplicity, we assume that $\mathbf{Set}_{\mathrm{A}_{1}}$ has at least two different assumption. The rest of the proof follows immediately.\par
To prove case (i), we consider now two distinct objects $\mathbf{Set}_{\mathrm{A}_{1}}$ and $\mathbf{Set}_{\mathrm{A}_{2}}$ given by $\mathbf{Set}_{\mathrm{A}_{1}}=\left\{\mathrm{A}_{1}^{(1)}, \mathrm{A}_{1}^{(2)}, \mathrm{A}_{1}^{(3)} \right\}$ and $\mathbf{Set}_{\mathrm{A}_{2}}=\left\{\mathrm{A}_{1}^{(1)}, \mathrm{A}_{1}^{(2)}, \mathrm{A}_{2}^{(1)} \right\}$, respectively. Similar to cases (ii) and (iii), we construct now two other objects in two different ways as follows:

\begin{equation*}
\begin{array}{lcl}
\displaystyle \mathbf{Set}_{\mathrm{A}_{3}} & = & \displaystyle \left\{\mathrm{A}_{1}^{(1)}, \mathrm{A}_{1}^{(2)}, \mathrm{A}_{1}^{(3)} \right\} \setminus \left\{\mathrm{A}_{1}^{(2)}, \mathrm{A}_{1}^{(3)} \right\} = \left\{\mathrm{A}_{1}^{(1)} \right\}, \\
\\
\displaystyle \mathbf{Set}_{\mathrm{A}_{4}} & = & \displaystyle \left\{\mathrm{A}_{1}^{(1)}, \mathrm{A}_{1}^{(2)}, \mathrm{A}_{1}^{(3)} \right\} \setminus \left\{\mathrm{A}_{1}^{(1)}, \mathrm{A}_{1}^{(3)} \right\} = \left\{\mathrm{A}_{1}^{(2)} \right\},
\end{array}
\end{equation*}
and 
\begin{equation*}
\begin{array}{lcl}
\displaystyle \mathbf{Set}_{\mathrm{A}_{3}} & = & \displaystyle \left\{\mathrm{A}_{1}^{(1)}, \mathrm{A}_{1}^{(2)}, \mathrm{A}_{2}^{(1)} \right\} \setminus \left\{\mathrm{A}_{1}^{(2)}, \mathrm{A}_{2}^{(1)} \right\} = \left\{\mathrm{A}_{1}^{(1)} \right\}, \\
\\
\displaystyle \mathbf{Set}_{\mathrm{A}_{4}} & = & \displaystyle \left\{\mathrm{A}_{1}^{(1)}, \mathrm{A}_{1}^{(2)}, \mathrm{A}_{2}^{(1)}\right\} \setminus \left\{\mathrm{A}_{1}^{(1)}, \mathrm{A}_{2}^{(1)} \right\} = \left\{ \mathrm{A}_{1}^{(2)} \right\}.
\end{array}
\end{equation*}
This construction is illustrated by the following diagram:

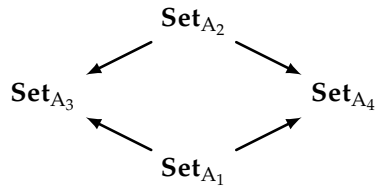
\begin{equation*}
\begin{tikzpicture}[scale=1,cap=round]
\tikzstyle{information text}=[rounded corners,fill=red!10,inner sep=1ex]
\begin{scope}[>=latex]
    \path (0,0) node (a) {$\mathbf{Set}_{\mathrm{A}_{1}}$};
    \path (0,2) node (b) {$\mathbf{Set}_{\mathrm{A}_{2}}$};
    \path (-2,1) node (c) {$\mathbf{Set}_{\mathrm{A}_{3}}$};
    \path (2,1) node (d) {$\mathbf{Set}_{\mathrm{A}_{4}}$};
    \draw[->,line width=1.0pt] (a) -- (c);
    \draw[->,line width=1.0pt] (a) -- (d); 
    \draw[->,line width=1.0pt] (b) -- (c); 
    \draw[->,line width=1.0pt] (b) -- (d); 
  \end{scope}
\end{tikzpicture}
\end{equation*}
Thus, the constructed category is partially ordered, and since objects $\mathbf{Set}_{\mathrm{A}_{1}}$ and $\mathbf{Set}_{\mathrm{A}_{2}}$ are not related, this category does not contain neither the most complex nor the simplest objects, since no object satisfies assumptions of Corollary~\ref{Corollary:simplest_most_complex}.\par
For proving case (iv), let us consider the object $\mathbf{Set}_{\mathrm{A}_{1}}=\left\{\mathrm{A}_{1}^{(1)}, \mathrm{A}_{1}^{(2)}, \mathrm{A}_{1}^{(3)}, \mathrm{A}_{1}^{(4)} \right\}$, and let us construct several other objects according to the following commutative diagram
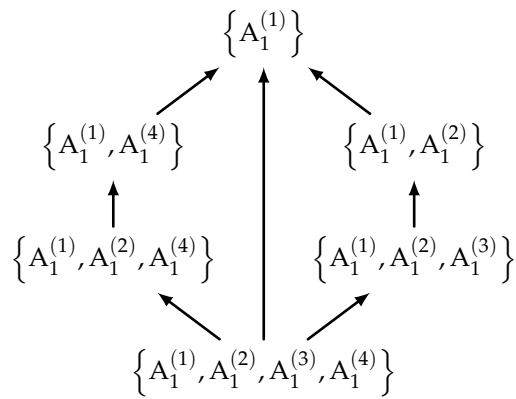
\begin{equation*}
\begin{tikzpicture}[scale=1,cap=round]
\tikzstyle{information text}=[rounded corners,fill=red!10,inner sep=1ex]
\begin{scope}[>=latex]
    \path (0,0) node (a) {$\left\{\mathrm{A}_{1}^{(1)}, \mathrm{A}_{1}^{(2)}, \mathrm{A}_{1}^{(3)}, \mathrm{A}_{1}^{(4)} \right\}$};
    \path (-2,1.5) node (b) {$\left\{\mathrm{A}_{1}^{(1)}, \mathrm{A}_{1}^{(2)}, \mathrm{A}_{1}^{(4)} \right\}$};
    \path (2,1.5) node (c) {$\left\{\mathrm{A}_{1}^{(1)}, \mathrm{A}_{1}^{(2)}, \mathrm{A}_{1}^{(3)} \right\}$};
    \path (-2,3) node (d) {$\left\{\mathrm{A}_{1}^{(1)}, \mathrm{A}_{1}^{(4)} \right\}$};
    \path (2,3) node (e) {$\left\{\mathrm{A}_{1}^{(1)}, \mathrm{A}_{1}^{(2)} \right\}$};
    \path (0,4.5) node (f) {$\left\{\mathrm{A}_{1}^{(1)}\right\}$};
    \draw[->,line width=1.0pt] (a) -- (b);
    \draw[->,line width=1.0pt] (a) -- (c); 
    \draw[->,line width=1.0pt] (b) -- (d); 
    \draw[->,line width=1.0pt] (d) -- (f); 
    \draw[->,line width=1.0pt] (c) -- (e); 
    \draw[->,line width=1.0pt] (e) -- (f); 
    \draw[->,line width=1.0pt] (a) -- (f); 
  \end{scope}
\end{tikzpicture}
\end{equation*}
Although the diagram is commutative, but the objects on the left side are not related to the objects of the right side in the sense of Definition~\ref{Defition:complexity}. Thus, we have a partially ordered category, where both the most complex object $\left\{\mathrm{A}_{1}^{(1)}\right\}$ and the simplest object $\left\{\mathrm{A}_{1}^{(1)}, \mathrm{A}_{1}^{(2)}, \mathrm{A}_{1}^{(3)}, \mathrm{A}_{1}^{(4)} \right\}$ exist simultaneously. Hence, the proposition is proved.\par
\end{proof}\par
Next, we have the following theorem:
\begin{Theorem}\label{Theorem:Partially_totally_ordered}
Consider a category $\mathbf{Model}_{1}$ with $n$ objects. If the most complex object $\mathbf{Set}_{\mathrm{A}_{1}}$ and the simplest object $\mathbf{Set}_{\mathrm{A}_{n}}$ exist simultaneously in the category $\mathbf{Model}_{1}$, then $\mathbf{Model}_{1}$ is either a totally ordered category, or contains at least two totally ordered subcategories.
\end{Theorem}
\begin{proof}
The proof of the theorem follows immediately from Corollary~\ref{Corollary:simplest_most_complex}, Proposition~\ref{Proposition:partially_simplest_complex}, and Definition~\ref{Definition:partially_totally_categories}. Looking at the proof of the case (iv) in Proposition~\ref{Proposition:partially_simplest_complex}, we see immediately that two totally ordered subcategories exist. The case of only one totally ordered subcategory is excluded by the assumption of simultaneous existence of the most complex and the simplest objects. Further, if the most complex and the simplest objects exist simultaneously and all objects in the category $\mathbf{Model}_{1}$ are related by help of complexity, then it follows immediately that $\mathbf{Model}_{1}$ is a totally ordered category.
\end{proof}\par
Evidently, the last statement can be straightforwardly generalised as follows:
\begin{Theorem}
Every partially ordered category of mathematical models contains at least one totally ordered category of mathematical models as a subcategory.
\end{Theorem}\par

\section{Convertible mathematical models}\label{Section:convertible}

In this section, we will discuss the mappings $S$ between sets of assumptions and the corresponding models appearing in Definition~\ref{Definition:category_of_mathematical_models}, and as we will see from the upcoming discussion, the role of mappings $S$ provides clear reasoning why objects of categories of mathematical models are sets of assumptions and not the models themselves. The mappings $S$ are generally not invertible, because they represent a formalisation process of basic modelling assumptions in terms of mathematical expressions. Moreover, these mappings are also not unique, since the same set of assumptions can be formalised differently. However, if objects in a category have been ordered (partially or totally) according their complexity, then the mappings will preserve this structure. Thus, these mappings are structure preserving mappings, i.e. they are functors.\par
Because the mappings between sets of assumptions and the corresponding mathematical models are functorial, then, in fact, the mathematical models constitute also a category. However, since final form of a model depends on the formalisation process, it is more difficult to work directly with categories of models, rather than to describe categories of sets of assumptions, as we have done already. Nonetheless, we will point out now some results related to the models directly. First, we summarise the above discussion in the following definition:
\begin{Definition}\label{Definition:convertible_models}
Let $\mathbf{Set}_{\mathrm{A}_{1}}$ be an object in the category $\mathbf{Model}_{1}$, and let $\mathrm{B}_{1}$ and $\mathrm{B}_{2}$ be two possible model formulations associated with the object $\mathbf{Set}_{\mathrm{A}_{1}}$ via two functors $F$ and $G$. Then the model formulations $\mathrm{B}_{1}$ and $\mathrm{B}_{2}$ are connected via a natural transformation of functors $\vartheta$, and the model formulations $\mathrm{B}_{1}$ and $\mathrm{B}_{2}$ are called {\bfseries convertible}. This construction corresponds to the commutative diagram

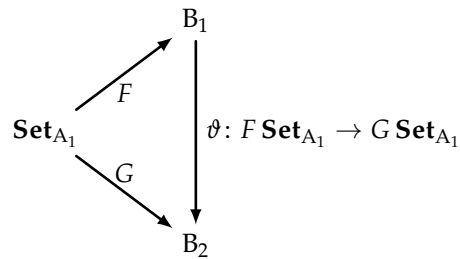
\begin{equation*}
\begin{tikzpicture}[scale=1,cap=round]
  \tikzstyle{information text}=[rounded corners,fill=red!10,inner sep=1ex]
  \begin{scope}[>=latex]
    \path (0,0) node (a) {$\mathbf{Set}_{\mathrm{A}_{1}}$};
    \path (2,1.5) node (b) {$\mathrm{B}_{1}$};
    \path (2,-1.5) node (c) {$\mathrm{B}_{2}$};
    \draw[->,line width=1.0pt] (a) -- (b) node[midway,below] {$F$};
    \draw[->,line width=1.0pt] (a) -- (c) node[midway,above] {$G$};
    \draw[->,line width=1.0pt] (b) -- (c) node[midway,right] {$\vartheta\colon F\, \mathbf{Set}_{\mathrm{A}_{1}} \to G\,  \mathbf{Set}_{\mathrm{A}_{1}}$};
  \end{scope}
\end{tikzpicture}
\end{equation*}
Moreover, models which are instantiated by convertible model formulations will be called {\bfseries convertible models}.
\end{Definition}
Obviously, because different model formulations are related to the same set of assumptions, the model complexity of these formulations remains the same. Thus, we have immediately the following corollary:
\begin{Corollary}
Convertible models have the same complexity.
\end{Corollary}\par
The discussion about convertible mathematical models underlines once more why sets of assumptions are considered as objects in categories of mathematical models, and not model formulations directly. Assume for a moment, that the latter would be the case and consider the following diagram with three objects for simplicity:

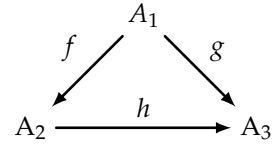
\begin{equation*}
\begin{tikzpicture}[scale=1,cap=round]
  \tikzstyle{information text}=[rounded corners,fill=red!10,inner sep=1ex]
  \begin{scope}[>=latex]
    \path (0,0) node (a) {${A}_{1}$};
    \path (-1.5,-1.5) node (b) {$\mathrm{A}_{2}$};
    \path (1.5,-1.5) node (c) {$\mathrm{A}_{3}$};
    \draw[->,line width=1.0pt] (a) -- (b) node[midway,above left] {$f$};
    \draw[->,line width=1.0pt] (a) -- (c) node[midway,above right] {$g$};
    \draw[->,line width=1.0pt] (b) -- (c) node[midway,above] {$h$};
  \end{scope}
\end{tikzpicture}
\end{equation*}
Moreover, assume additionally that the model formulations $A_{1}$ and $A_{2}$ are convertible in the sense of Definition~\ref{Definition:convertible_models}, while the model formulation $A_{3}$ is not associated with the same set of assumptions. Thus, we would end up with two kinds of morphisms in the category: morphism $f$ plays the same role as the natural transformation $\vartheta$ in Definition~\ref{Definition:convertible_models}, while morphisms $g$ and $h$ represent complexity-relation on the level of model formulations. Obviously, it is necessary to be able to distinguish between the two kinds of morphisms, which would imply much more complicated constructions for the structure of the category, as well as for relations between its objects.\par
As a simple immediate example indicating the necessity for considering convertible mathematical models, let us consider the classical model of linear elasticity describing deformations of an elastic body in a static case. The classical formulation of this model is given by the following system of equations
\begin{equation}
\label{elasticity_tensor_version}
\left\{\begin{array}{l}
\displaystyle \mathrm{div}\,\tilde{\boldsymbol\sigma}+\rho\,\mathbf{K}=0, \\
\displaystyle \tilde{\boldsymbol\varepsilon}=\frac{1}{2}\left[\nabla\mathbf{u}+\left(\nabla\mathbf{u}\right)^{\mathrm{T}}\right], \\
\displaystyle \tilde{\boldsymbol\sigma}=2\mu\left(\frac{\nu}{1-2\nu}\vartheta\,\tilde{\mathbf{E}} + \tilde{\boldsymbol\varepsilon}\right),
\end{array}\right. \displaystyle \vartheta=\mathrm{div}\,\mathbf{u} = \frac{\partial u_{1}}{\partial x_{1}} + \frac{\partial u_{2}}{\partial x_{2}} + \frac{\partial u_{3}}{\partial x_{3}},
\end{equation}
where $\tilde{\boldsymbol\sigma}$ is a symmetric stress tensor, $\tilde{\boldsymbol\varepsilon}$ is a symmetric strain tensor, $\mathbf{u}$ is a displacement vector, $\rho$ is a material density, $\nu$ is the Poisson's ration, and $\mathbf{K}$ is the volume force. System of equations~(\ref{elasticity_tensor_version}) is the classical tensor version of elasticity equations, see for example \cite{Lurie}. However, the Lam\'{e} equation
\begin{equation}
\label{lame_equation}
\mu\,\Delta\mathbf{u}+(\lambda+\mu)\mathrm{grad}\,\,\mathrm{div}\,\mathbf{u} +\rho\mathbf{K}=0,
\end{equation}
is often used in practice as well. Furthermore, model of linear elasticity can be also written as follows
\begin{equation}
\label{elasticity_quaternionic}
D M D u = 0, \mbox{ with } D = \sum\limits_{k=1}^{3}\mathbf{e}_{k}\partial_{k}, \mbox{ and } u=u_{0}+\mathbf{u},
\end{equation}
where the multiplication operator $M$ is defined by
\begin{equation*}
M\mathbf{u} := \frac{m-2}{2(m-1)}u_{0}+\mathbf{u}, \quad m:=\nu^{-1}.
\end{equation*}
Equation~(\ref{elasticity_quaternionic}) is a quaternionic form of elasticity model with $D$ denoting the Dirac operator, see \cite{Guerlebeck_2} for all details on quaternionic analysis and its applications.\par
For the sake of clarity of further considerations, let us denote the models~(\ref{elasticity_tensor_version})-(\ref{elasticity_quaternionic}) as follows:
\begin{equation*}
\begin{array}{ll}
\displaystyle\mathrm{B}_{1} := & \displaystyle \left\{\begin{array}{l}
\displaystyle \mathrm{div}\,\tilde{\boldsymbol\sigma}+\rho\,\mathbf{K}=0, \\
\displaystyle \tilde{\boldsymbol\varepsilon}=\frac{1}{2}\left[\nabla\mathbf{u}+\left(\nabla\mathbf{u}\right)^{\mathrm{T}}\right], \\
\displaystyle \tilde{\boldsymbol\sigma}=2\mu\left(\frac{\nu}{1-2\nu}\vartheta\,\tilde{\mathbf{E}} + \tilde{\boldsymbol\varepsilon}\right),
\end{array}\right. \displaystyle \vartheta=\mathrm{div}\,\mathbf{u} = \frac{\partial u_{1}}{\partial x_{1}} + \frac{\partial u_{2}}{\partial x_{2}} + \frac{\partial u_{3}}{\partial x_{3}}, \\
\displaystyle\mathrm{B}_{2} := & \displaystyle \mu\,\Delta\mathbf{u}+(\lambda+\mu)\mathrm{grad}\,\,\mathrm{div}\,\mathbf{u} +\rho\mathbf{K}=0, \\
\displaystyle\mathrm{B}_{3} := & \displaystyle D M D u = 0, \mbox{ with } D = \sum\limits_{k=1}^{3}\mathbf{e}_{k}\partial_{k}, \mbox{ and } u=u_{0}+\mathbf{u}.
\end{array}
\end{equation*}
A possible representation of these models is provided by the diagram
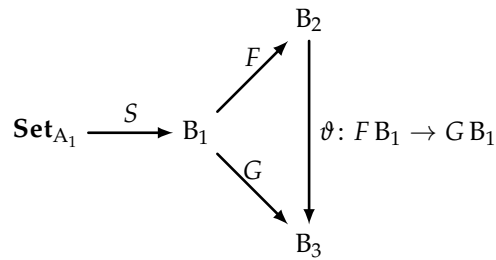
\begin{equation*}
\begin{tikzpicture}[scale=1,cap=round]
  \tikzstyle{information text}=[rounded corners,fill=red!10,inner sep=1ex]
  \begin{scope}[>=latex]
    \path (-2,0) node (a) {$\mathbf{Set}_{\mathrm{A}_{1}}$};
    \path (0,0) node (b) {$\mathrm{B}_{1}$};
    \path (1.5,1.5) node (c) {$\mathrm{B}_{2}$};
    \path (1.5,-1.5) node (d) {$\mathrm{B}_{3}$};
    \draw[->,line width=1.0pt] (a) -- (b) node[midway,above] {$S$};
    \draw[->,line width=1.0pt] (b) -- (c) node[midway,above] {$F$};
    \draw[->,line width=1.0pt] (b) -- (d) node[midway,above] {$G$};
    \draw[->,line width=1.0pt] (c) -- (d) node[midway,right] {$\vartheta\colon F\, \mathrm{B}_{1} \to G\, \mathrm{B}_{1}$};
  \end{scope}
\end{tikzpicture}
\end{equation*}
Here, functor $S$ is a formalisation process of basic set of assumptions of linear elasticity $\mathbf{Set}_{\mathrm{A}_{1}}$ in the tensor form of model formulation $\mathrm{B}_{1}$, after that, the tensor form can be further reformulated into the Lam\'e equation $\mathrm{B}_{2}$, or into the quaternionic form $\mathrm{B}_{3}$ via functorial mappings $F$ and $G$. In some sense, the above diagram reflects traditional way of developing different model formulations: at first, the original form is introduced, and after that, several more specific forms better suitable for selected methods are introduced. Moreover, looking in particular at the quaternionic formulation $B_{3}$, it becomes clear that this form is not obtained directly through the formalisation process of $\mathbf{Set}_{\mathrm{A}_{1}}$ (at least no quaterninic-based modelling of linear elasticity has been reported till now), but through reformulation of either Lam\'e equation or the tensor form, see again \cite{Guerlebeck_2}.\par

\section{Illustrative examples}\label{Section:example}

In this section, we illustrate the constructions of category theory-based modelling methodology presented in previous sections on two examples: first, we discuss classical models of beam theories, and after that, we discuss aerodynamic models used in bridge engineering. These examples have been already presented in works \cite{Guerlebeck,Kavrakov} at the time of first steps towards developing the category theory-based modelling methodology. Therefore, it is necessary to revisit these examples for underlying further development of the theory.\par

\subsection{Categorical modelling of beam theories}

Transverse vibrations of one-dimensional beams are typically modelled by one of three common beam theories: Bernoulli-Euler theory, Rayleigh theory, and Timoshenko theory. Thus, let us consider a category of mathematical models, denoted by $\mathbf{Beam}$, containing as objects sets of assumptions $\mathbf{Set}_{\mathrm{B-E}}$, $\mathbf{Set}_{\mathrm{R}}$, $\mathbf{Set}_{\mathrm{T}}$ corresponding to the Bernoulli-Euler, Rayleigh, and Timoshenko beam theories, respectively. We start our discussion on the construction of category $\mathbf{Beam}$ by explicitly listing the sets of assumptions, which are given in Table \ref{table_assumptions_beams}.
\begin{specialtable}[H] 
\caption{Sets of assumptions of beam theories\label{table_assumptions_beams}}
\begin{tabular}{|p{9cm}|c|c|c|}
\hline
\centering {\bfseries Assumptions} & $\mathbf{Set}_{\mathrm{B-E}}$ & $\mathbf{Set}_{\mathrm{R}}$ & $\mathbf{Set}_{\mathrm{T}}$ \\
\hline
1.~Cross sections of a beam that are planes remain planes after the deformation process & $\boldsymbol{+}$ & $\boldsymbol{+}$ & $\boldsymbol{+}$ \\
2.~Normal stresses on planes parallel to the axis of a beam are infinitesimal & $\boldsymbol{+}$ & $\boldsymbol{+}$ & $\boldsymbol{+}$ \\
3.~A beam has a constant cross section & $\boldsymbol{+}$ & $\boldsymbol{+}$ & $\boldsymbol{+}$ \\
4.~A beam is made of a homogeneous isotropic material & $\boldsymbol{+}$ & $\boldsymbol{+}$ & $\boldsymbol{+}$ \\
5.~Cross sections of a beam perpendicular to its axis remain perpendicular to the deformed axis & $\boldsymbol{+}$ & $\boldsymbol{+}$ &  \\
6.~Rotation inertia of cross sections of a beam is omitted & $\boldsymbol{+}$ &  &  \\
\hline
\end{tabular}
\end{specialtable}\par
\begin{Remark}
The assumptions, as listed in Table~\ref{table_assumptions_beams}, are formulated by help of natural language, however in some cases it is more convenient to formulate sets of assumptions directly in terms of mathematical expressions, or as a mixture of both. Although from the set-theoretic point of view such a freedom in writing sets of assumptions is not completely justified, it is acceptable in our setting because each set of assumption written in natural language can be rigorously formalised in terms of mathematical expressions. Thus, writing mathematical expressions in sets of assumptions can be considered as a kind of {\bfseries syntactic sugar}, similar to programming languages terminology. Of course, this analogy not perfect but reflects a general point of view on writing sets of assumptions.
\end{Remark}\par
Since derivation of beam models is well known, it will be omitted. Set of assumption $\mathbf{Set}_{\mathrm{B-E}}$ of the Bernoulli-Euler theory leads to the following beam equation:
\begin{equation*}
\rho\,F\frac{\partial^{2}u}{\partial t^{2}} + E\,I_{y}\frac{\partial^{4} u}{\partial x^{4}} = 0,
\end{equation*}
where $E$ is the Young's modulus of the material, $I_{y}$ is the moment of inertia, $rho$ is the density of material, and $F$ is the area of cross section. Next, set of assumption $\mathbf{Set}_{\mathrm{R}}$ of the Rayleigh theory leads to the equation:
\begin{equation*}
\rho\,F\frac{\partial^{2}u}{\partial t^{2}} + E\,I_{y}\frac{\partial^{4} u}{\partial x^{4}} - \rho\,I_{y}\frac{\partial^{4} u}{\partial x^{2}\partial t^{2}} = 0.
\end{equation*}
Finally, if the effect of bending of cross sections is taken into account, then set of assumption $\mathbf{Set}_{\mathrm{T}}$ of the Timoshenko theory is obtained, which leads to the system of differential equations:
\begin{equation*}
\left\{\begin{array}{ccl}
\displaystyle \rho F \frac{\partial^{2} u}{\partial t^{2}} - \aleph\mu F \frac{\partial^{2} u}{\partial x^{2}} + \aleph\mu F \frac{\partial \varphi}{\partial x} & = & 0, \\
\\
\displaystyle \rho I_{y} \frac{\partial^{2} \varphi}{\partial t^{2}} - E I_{y} \frac{\partial^{2} \varphi}{\partial x^{2}} + \aleph\mu F\left( \varphi - \frac{\partial u}{\partial x}\right) & = & 0,
\end{array}\right.
\end{equation*}
where $\varphi$ is the angle of rotation of the normal to the mid-surface of the beam, $\aleph$ is the Timoshenko shear coefficient, which depends on the geometry of the beam, and $\mu$ is the shear modulus. After some calculations this system can be reformulated in terms of only one partial differential equation for $u$ as follows:
\begin{equation*}
\rho\,F\frac{\partial^{2}u}{\partial t^{2}} + E\,I_{y}\frac{\partial^{4} u}{\partial x^{4}} - \rho\,I_{y}\left(1 + \frac{E}{\aleph\,\mu} \right)\frac{\partial^{4} u}{\partial x^{2}\partial t^{2}} + \frac{\rho^{2}I_{y}}{\aleph\,\mu}\frac{\partial^{4} u}{\partial t^{4}} = 0.
\end{equation*}\par
Looking at the above beam models from the categorical perspective, we can summarise these models and their sets of assumptions as follows:
\begin{equation*}
\begin{array}{ccll}
\displaystyle \mathbf{Set}_{\mathrm{B-E}} & \displaystyle \stackrel{S}{\mapsto} & \displaystyle \rho\,F\frac{\partial^{2}u}{\partial t^{2}} + E\,I_{y}\frac{\partial^{4} u}{\partial x^{4}} = 0 & \displaystyle =: \mathbf{A}, \\
\\
\displaystyle \mathbf{Set}_{\mathrm{R}} & \displaystyle \stackrel{S}{\mapsto} & \displaystyle \rho\,F\frac{\partial^{2}u}{\partial t^{2}} + E\,I_{y}\frac{\partial^{4} u}{\partial x^{4}} - \rho\,I_{y}\frac{\partial^{4} u}{\partial x^{2}\partial t^{2}} = 0 & \displaystyle =: \mathbf{B},\\
\\
\displaystyle \mathbf{Set}_{\mathrm{T}} & \displaystyle \stackrel{S}{\mapsto} & \displaystyle \rho\,F\frac{\partial^{2}u}{\partial t^{2}} + E\,I_{y}\frac{\partial^{4} u}{\partial x^{4}} - \rho\,I_{y}\left(1 + \frac{E}{\aleph\,\mu} \right)\frac{\partial^{4} u}{\partial x^{2}\partial t^{2}} + \frac{\rho^{2}I_{y}}{\aleph\,\mu}\frac{\partial^{4} u}{\partial t^{4}} = 0 & \displaystyle =: \mathbf{C}_{1}, \\
\\
\displaystyle \mathbf{Set}_{\mathrm{T}} & \displaystyle \stackrel{S}{\mapsto} & \left\{\begin{array}{ccl}
\displaystyle \rho F \frac{\partial^{2} u}{\partial t^{2}} - \aleph\mu F \frac{\partial^{2} u}{\partial x^{2}} + \aleph\mu F \frac{\partial \varphi}{\partial x} & = & 0, \\
\\
\displaystyle \rho I_{y} \frac{\partial^{2} \varphi}{\partial t^{2}} - E I_{y} \frac{\partial^{2} \varphi}{\partial x^{2}} + \aleph\mu F\left( \varphi - \frac{\partial u}{\partial x}\right) & = & 0.
\end{array}\right. & \displaystyle =: \mathbf{C}_{2}, 
\end{array}
\end{equation*}
where $S$ are formalisation mappings, as discussed before. It is worth making the remark:
\begin{Remark}
Note that, in general, mappings $S$ can be different for each set of assumptions, or, can be the same if all equations are derived based on the same principle, e.g. the Hamilton's principle. If the fact that different formalisation processes have been used to obtain models from the sets of assumptions in one category is essential for the application, then it is necessary to indicate this fact by using sub-scripts, i.e. $S_{1}$, $S_{2}$, $\ldots$, otherwise the general notation for the formalisation mappings might be kept.
\end{Remark}\par
By using Definition~\ref{Defition:complexity}, the category $\mathbf{Beam}$ can be straightforwardly equipped with the commutative diagram
\begin{equation*}
\begin{tikzpicture}[scale=1,cap=round]
  \tikzstyle{information text}=[rounded corners,fill=red!10,inner sep=1ex]
  \begin{scope}[>=latex]
    \path (0,0) node (a) {$\mathbf{Set}_{\mathrm{B-E}}$};
    \path (2,1.5) node (b) {$\mathbf{Set}_{\mathrm{R}}$};
    \path (4,0) node (c) {$\mathbf{Set}_{\mathrm{T}}$};
    \draw[->,line width=1.0pt] (a) -- (b) node[midway,above left] {$f$};
    \draw[->,line width=1.0pt] (b) -- (c) node[midway,above right] {$g$};
    \draw[->,line width=1.0pt] (a) -- (c) node[midway,below] {$h=g\circ f$};
  \end{scope}
\end{tikzpicture}
\end{equation*}
The morphisms $f$, $g$, and $h$ indicate the simple fact, that one beam theory can be obtained from another by weakening basic assumptions. Moreover, the above diagram clearly indicate that the object $\mathbf{Set}_{\mathrm{T}}$ (Timoshenko theory) is the most complex, the object $\mathbf{Set}_{\mathrm{R}}$ (Rayleigh theory) has higher complexity than the object $\mathbf{Set}_{\mathrm{B-E}}$ (Bernoulli-Euler theory), which is the simplest object. The same ordering holds for the corresponding model instantiations. Next, let us list the following facts we know about the category $\mathbf{Beam}$:
\begin{itemize}
\item it is a totally ordered category;
\item the object $\mathbf{Set}_{\mathrm{B-E}}$ is the initial object of this category;
\item the object $\mathbf{Set}_{\mathrm{T}}$ is the terminal object of this category;
\item models $\mathbf{C}_{1}$ and $\mathbf{C}_{2}$ are convertible, since they represent different formulations of the assumptions of Timoshenko theory.
\end{itemize}
Note that, first three facts, as well as the commutative diagram presented above, do not require, in fact, models themself, because these facts are solely obtained simply from the sets of assumptions, i.e. by looking at the objects in category $\mathbf{Beam}$. Thus, the categorical point of view introduced in the previous section reflects the following idea:\\[7pt]
{\itshape The principle difference between models lies not in their final form, but in the basic modelling assumptions these models constructed from.}\\[7pt] 
Finally, let us look at the level of models, where the following diagram is obtained
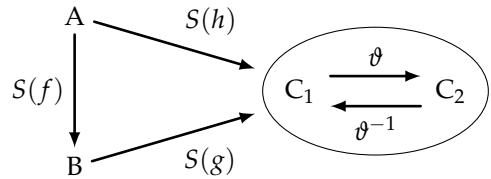
\begin{equation*}
\begin{tikzpicture}[scale=1,cap=round]
  \tikzstyle{information text}=[rounded corners,fill=red!10,inner sep=1ex]
  \begin{scope}[>=latex]
    \path (0,0) node (a) {$\mathrm{A}$};
    \path (0,-2) node (b) {$\mathrm{B}$};
    \path (3,-1.0) node (c) {$\mathrm{C}_{1}$};
    \path (5,-1.0) node (d) {$\mathrm{C}_{2}$};
    \draw (4,-1.0) circle [x radius=1.5cm, y radius=0.85cm];
    \draw[->,line width=1.0pt] (a) -- (b) node[midway,left] {$S(f)$};
    \draw[->,line width=1.0pt] (3.4,-0.8) -- (4.6,-0.8) node[midway,above] {$\vartheta$};
    \draw[<-,line width=1.0pt] (3.4,-1.2) -- (4.6,-1.2) node[midway,below] {$\vartheta^{-1}$};
    \draw[->,line width=1.0pt] (b) -- (2.4,-1.3) node[midway,below right] {$S(g)$};
    \draw[->,line width=1.0pt] (a) -- (2.4,-0.7) node[midway,above right] {$S(h)$};
  \end{scope}
\end{tikzpicture}
\end{equation*}
where $\vartheta$ denotes a natural transformation appearing in the definition of convertible models, recall Definition~\ref{Definition:convertible_models}.\par

\subsection{Category of aerodynamic models revisited}

Next, we briefly revisit the example of aerodynamic models used in bridge engineering presented in \cite{Kavrakov}. Since the idea is only briefly discuss categorical constructions introduced in previous sections, we will not present aerodynamic models in details, but we refer to works \cite{Kavrakov_1,Kavrakov_2}. We consider the category $\mathbf{AeroModel}$ containing as objects the following sets of assumptions of mathematical models: (i) ${\mathbf{ST}}$ (steady model); (ii) ${\mathbf{LST}}$ (linear steady model); (iii) ${\mathbf{QS}}$ (quasi-steady model); (iv) ${\mathbf{LQS}}$ (linear quasi-steady model); (v) ${\mathbf{LU}}$ (linear unsteady model); (vi) ${\mathbf{MQS}}$ (modified quasi-steady model); (vii) ${\mathbf{MBM}}$ (mode-by-mode model); (viii) ${\mathbf{CQS}}$ (corrected quasi-steady model); (ix) ${\mathbf{HNL}}$ (hybrid nonlinear model); (x) ${\mathbf{MNL}}$ (modified nonlinear model); and, (xi) ${\mathbf{NLU}}$ (nonlinear unsteady model). The structure of category $\mathbf{AeroModel}$ is provided by the following diagram (adapted from \cite{Kavrakov}):

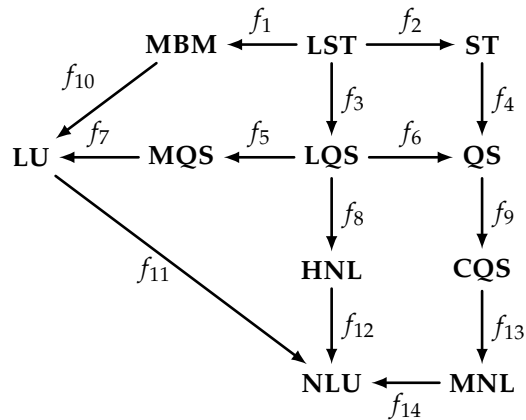
\begin{equation*}
\centering
\begin{tikzpicture}[scale=1,cap=round]
\tikzstyle{information text}=[rounded corners,fill=red!10,inner sep=1ex]
\begin{scope}[>=latex]  
\path (0,3) node (a) {${\mathbf{LST}}$};       
\path (2,3) node (b) {${\mathbf{ST}}$};        
\path (2,1.5) node (c) {${\mathbf{QS}}$};
\path (2,0) node (d) {${\mathbf{CQS}}$};   
\path (0,1.5) node (e) {${\mathbf{LQS}}$};       
\path (0,0) node (f) {${\mathbf{HNL}}$};   
\path (-2,1.5) node (g) {${\mathbf{MQS}}$};
\path (0,-1.5) node (h) {${\mathbf{NLU}}$};
\path (-2,3) node (i) {${\mathbf{MBM}}$};      
\path (-4,1.5) node (j) {${\mathbf{LU}}$}; 
\path (2,-1.5) node (m) {${\mathbf{MNL}}$};
\draw[->,line width=1.0pt] (a) -- (e) node[midway,right] {$f_{3}$}; 
\draw[->,line width=1.0pt] (a) -- (b) node[midway,above] {$f_{2}$};
\draw[->,line width=1.0pt] (b) -- (c) node[midway,right] {$f_{4}$};   
\draw[->,line width=1.0pt] (c) -- (d) node[midway,right] {$f_{9}$};    
\draw[->,line width=1.0pt] (e) -- (f) node[midway,right] {$f_{8}$};
\draw[->,line width=1.0pt] (e) -- (c) node[midway,above] {$f_{6}$};
\draw[->,line width=1.0pt] (f) -- (h) node[midway,right] {$f_{12}$};
\draw[->,line width=1.0pt] (d) -- (m) node[midway,right] {$f_{13}$};
\draw[->,line width=1.0pt] (j) -- (h) node[midway,left]  {$f_{11}$};
%\draw[->,line width=1.0pt] (e) -- (j) node[midway,above] {$f_{7}$};   
\draw[->,line width=1.0pt] (g) -- (j) node[midway,above] {$f_{7}$}; 
\draw[->,line width=1.0pt] (e) -- (g) node[midway,above] {$f_{5}$};  
\draw[->,line width=1.0pt] (a) -- (i) node[midway,above] {$f_{1}$};  
\draw[->,line width=1.0pt] (i) -- (j) node[midway,above left] {$f_{10}$};  
\draw[->,line width=1.0pt] (m) -- (h) node[midway,below] {$f_{14}$};                         
\end{scope}
\end{tikzpicture}
\end{equation*}

Let us now list some facts we know about the category $\mathbf{AeroModel}$:
\begin{itemize}
\item it is a partially ordered category;
\item the object ${\mathbf{LST}}$ is the initial object of this category;
\item the object ${\mathbf{NLU}}$ is the terminal object of this category;
\item according to Theorem~\ref{Theorem:Partially_totally_ordered} several totally ordered subcategories exists, which are

\begin{equation*}
\centering
\begin{tikzpicture}[scale=1,cap=round]
\tikzstyle{information text}=[rounded corners,fill=red!10,inner sep=1ex]
\begin{scope}[>=latex]  
\path (-0.6,0) node {$1.$};       
\path (0,0) node (a1) {${\mathbf{LST}}$};       
\path (2,0) node (b1) {${\mathbf{ST}}$};        
\path (4,0) node (c1) {${\mathbf{QS}}$};
\path (6,0) node (d1) {${\mathbf{CQS}}$};   
\path (8,0) node (e1) {${\mathbf{MNL}}$};       
\path (10,0) node (f1) {${\mathbf{NLU}}$};   
\draw[->,line width=1.0pt] (a1) -- (b1) node[midway,above] {$f_{2}$};
\draw[->,line width=1.0pt] (b1) -- (c1) node[midway,above] {$f_{4}$};   
\draw[->,line width=1.0pt] (c1) -- (d1) node[midway,above] {$f_{9}$};    
\draw[->,line width=1.0pt] (d1) -- (e1) node[midway,above] {$f_{13}$};
\draw[->,line width=1.0pt] (e1) -- (f1) node[midway,above] {$f_{14}$};

\path (-0.6,-1) node {$2.$};       
\path (0,-1) node (a2) {${\mathbf{LST}}$};       
\path (2,-1) node (b2) {${\mathbf{LQS}}$};        
\path (4,-1) node (c2) {${\mathbf{QS}}$};
\path (6,-1) node (d2) {${\mathbf{CQS}}$};   
\path (8,-1) node (e2) {${\mathbf{MNL}}$};       
\path (10,-1) node (f2) {${\mathbf{NLU}}$};   
\draw[->,line width=1.0pt] (a2) -- (b2) node[midway,above] {$f_{3}$};
\draw[->,line width=1.0pt] (b2) -- (c2) node[midway,above] {$f_{6}$};   
\draw[->,line width=1.0pt] (c2) -- (d2) node[midway,above] {$f_{9}$};    
\draw[->,line width=1.0pt] (d2) -- (e2) node[midway,above] {$f_{13}$};
\draw[->,line width=1.0pt] (e2) -- (f2) node[midway,above] {$f_{14}$};

\path (-0.6,-2) node {$3.$};       
\path (0,-2) node (a3) {${\mathbf{LST}}$};       
\path (2,-2) node (b3) {${\mathbf{LQS}}$};        
\path (4,-2) node (c3) {${\mathbf{HNL}}$};
\path (6,-2) node (d3) {${\mathbf{NLU}}$};   
\draw[->,line width=1.0pt] (a3) -- (b3) node[midway,above] {$f_{3}$};
\draw[->,line width=1.0pt] (b3) -- (c3) node[midway,above] {$f_{8}$};   
\draw[->,line width=1.0pt] (c3) -- (d3) node[midway,above] {$f_{12}$};    

\path (-0.6,-3) node {$4.$};       
\path (0,-3) node (a4) {${\mathbf{LST}}$};       
\path (2,-3) node (b4) {${\mathbf{LQS}}$};        
\path (4,-3) node (c4) {${\mathbf{MQS}}$};
\path (6,-3) node (d4) {${\mathbf{LU}}$};   
\path (8,-3) node (e4) {${\mathbf{NLU}}$};       
\draw[->,line width=1.0pt] (a4) -- (b4) node[midway,above] {$f_{3}$};
\draw[->,line width=1.0pt] (b4) -- (c4) node[midway,above] {$f_{5}$};   
\draw[->,line width=1.0pt] (c4) -- (d4) node[midway,above] {$f_{7}$};    
\draw[->,line width=1.0pt] (d4) -- (e4) node[midway,above] {$f_{11}$};

\path (-0.6,-4) node {$5.$};       
\path (0,-4) node (a5) {${\mathbf{LST}}$};       
\path (2,-4) node (b5) {${\mathbf{MBM}}$};        
\path (4,-4) node (c5) {${\mathbf{LU}}$};
\path (6,-4) node (d5) {${\mathbf{NLU}}$};   
\draw[->,line width=1.0pt] (a5) -- (b5) node[midway,above] {$f_{1}$};
\draw[->,line width=1.0pt] (b5) -- (c5) node[midway,above] {$f_{10}$};   
\draw[->,line width=1.0pt] (c5) -- (d5) node[midway,above] {$f_{11}$};    
\end{scope}
\end{tikzpicture}
\end{equation*}

\end{itemize}
Additionally, we can say that no models associated to the objects of $\mathbf{AeroModel}$ are convertible, but for that it is necessary to take a look at the derivation of models, see again \cite{Kavrakov} and references therein.\par

\section{Further characterisations of mathematical models and conclusions}\label{Section:conclusions}

In this section, we present some further ideas on characterisations of mathematical models. One of the most important aspect of applications of category theory is a definition of a {\itshape universal mapping property} (UMP), or simply, a {\itshape universal arrow}, which provides, in fact, a categorical characterisation of objects, see \cite{Awodey,MacLane} for details. Hence, it is important to discuss the universal arrow definition also in the context of category theory-based modelling methodology.\par
Let us consider a formalisation functor $S\colon \mathbf{Model} \to \mathbf{M}$, where $\mathbf{M}$ denotes formally a category of instantiations of mathematical models corresponding to the objects in $\mathbf{Model}$. Let $m$ be an object of $\mathbf{M}$, then a universal arrow from $m$ to $S$ is a pair $\langle r,u \rangle$ consisting of an object $r$ of $\mathbf{Model}$ and an arrow $r\colon m \to Sr$ of $\mathbf{M}$, such that to every pair $\langle d,f \rangle$ with $d$ an object of $\mathbf{Model}$ and $f\colon c \to Sd$ an arrow of $\mathbf{M}$, there is a unique arrow $f'\colon r \to d$ of $\mathbf{Model}$ with $Sf'\circ u = f$. Practical meaning of a universal arrow in the context of category theory-based modelling methodology is that to the same set of assumption can correspond only convertible model formulations.\par
Finally, we would like to provide another possible definition of a mathematical model in general, which would summarise our discussion in this paper:
\begin{Definition}\label{Definition:model_triple}
A mathematical model $\mathfrak{M}$ is a triple $\mathfrak{M}=\langle \mathbf{Set}, \mathcal{M}, S \rangle$, where
\begin{itemize}
\item $\mathbf{Set}$ is the set of assumptions of the model;
\item $\mathcal{M}$ is an instantiation of the model in terms of mathematical expressions and equations;
\item $S$ is a formalisation mapping, which formalises the set of assumptions $\mathbf{Set}$ into the model instantiation $\mathcal{M}$.
\end{itemize}
\end{Definition}
Relations between the models can be introduced again by help of Definition~\ref{Defition:complexity}. Definition~\ref{Definition:model_triple} proposes an abstract description of a mathematical model similar to the abstract algebraic approach presented in \cite{Legatiuk_2}. Thus, a connection between the category theory-based modelling methodology and abstract algebraic approach is established. Hence, both approaches to the modelling process in engineering might complement each other, and therefore, the connection between both approaches will be studied in future research.\par
In this paper, we have revisited the category theory-based modelling methodology proposed in recent years. The main idea of this modelling methodology is representation of mathematical models by help of categorical constructions. We have presented revised results from previous works, as well as new results and ideas supporting a deeper understanding of the modelling process in engineering. Moreover, two illustrative practical examples, namely categorical perspective of beam models and on aerodynamic models from bridge engineering, have been revisited. As it can be clearly seen from the examples, the category theory-based modelling methodology presented in this paper is indeed applicable in practice and provides various characterisations of mathematical models, relations between them, and final formulations of models. Finally, we have describe a universal arrow in the framework of the proposed modelling methodology.\par
Additionally, we would like to remark how the category theory-based modelling methodology presented in this paper can be used in a model selection process. After constructing a category of mathematical models, we can formulate criteria which must be satisfied by a model for a given practical problem, and thus, a subcategory of models satisfying these criteria can be constructed. Because we are on the abstract level of models, it is difficult to introduce a quantifiable criterion for the optimal model choice. Nonetheless, on the abstract level, the simplest model satisfying the criteria can be regarded as \textquotedblleft the optimal choice\textquotedblright\, in this case, because generally there is no need for overcomplicating the model. Furthermore, the difference in model assumptions, and thus, in model complexity, can be quantified by help of numerical calculations, as it has been illustrated in \cite{Kavrakov} for the case of aerodynamic models.\par
The scope of future research is related to a revision and deeper understanding of coupled mathematical models. A categorical description of coupled mathematical model will be using constructions and ideas introduced in this paper. However, due to more complex nature of coupled models, it is expected that more refined and advanced constructions will be necessary for a proper description of such models. Moreover, further ideas on a formal model comparison and model selection procedure, as well as a more strict approach to formulation of sets of assumptions, will be considered in future work.\par

%%%%%%%%%%%%%%%%%%%%%%%%%%%%%%%%%%%%%%%%%%
\vspace{6pt} 

%%%%%%%%%%%%%%%%%%%%%%%%%%%%%%%%%%%%%%%%%%
\funding{This research is supported by the German Research Foundation (DFG) through grant LE 3955/4-1.}

\acknowledgments{I would like to thank the reviewers for very helpful comments, which help not only improving the paper, but also brought new ideas for future research.}

%%%%%%%%%%%%%%%%%%%%%%%%%%%%%%%%%%%%%%%%%%
%% Optional
%\abbreviations{Abbreviations}{
%The following abbreviations are used in this manuscript:\\
%
%\noindent 
%\begin{tabular}{@{}ll}
%MDPI & Multidisciplinary Digital Publishing Institute\\
%DOAJ & Directory of open access journals\\
%TLA & Three letter acronym\\
%LD & Linear dichroism
%\end{tabular}}

%%%%%%%%%%%%%%%%%%%%%%%%%%%%%%%%%%%%%%%%%%
%% Optional
\appendixtitles{yes} % Leave argument "no" if all appendix headings stay EMPTY (then no dot is printed after "Appendix A"). If the appendix sections contain a heading then change the argument to "yes".
\appendixstart
\appendix

\section{Some basic definitions from category theory}

Following the classical works in category theory \cite{Awodey,MacLane}, we list here few important definitions.
\begin{Definition}
A {\itshape category} consists of the following data:
\begin{itemize}
\item {\itshape Objects}: $A,B,C,\ldots$
\item {\itshape Arrows (morphisms)}: $f,g,h,\ldots$
\item For each arrow $f$, there are given objects $\mathrm{dom}(f)$ and $\mathrm{cod}(f)$ called the {\itshape domain} and {\itshape codomain} of $f$, respectively. We write
\begin{equation*}
f\colon A\longrightarrow B \quad \mbox{or} \quad A \stackrel{f}{\longrightarrow} B
\end{equation*}
to indicate that $A=\mathrm{dom}(f)$ and $B=\mathrm{cod}(f)$.
\item Given arrows $f\colon A\longrightarrow B$ and $g\colon B\longrightarrow C$, that is, with $\mathrm{cod}(f)=\mathrm{dom}(g)$, there is given an arrow
\begin{equation*}
g\circ f\colon A\longrightarrow C
\end{equation*}
called the {\itshape composite} of $f$ and $g$.
\item For each object $A$, there is given an arrow
\begin{equation*}
1_{A}\colon A\longrightarrow A
\end{equation*}
called the {\itshape identity arrow} of $A$.
\end{itemize}
These data are required to satisfy the following laws:
\begin{itemize}
\item Associativity: $h\circ\left(g\circ f\right) = \left(h\circ g\right)\circ f$ for all $f\colon A\longrightarrow B$, $g\colon B\longrightarrow C$, $h\colon C\longrightarrow D$.
\item Unit: $f\circ 1_{A} = f = 1_{B}\circ f$ for all $f\colon A\longrightarrow B$.
\end{itemize}
\end{Definition}
\begin{Definition}
A {\itshape functor}
\begin{equation*}
F\colon \mathbf{C}\longrightarrow\mathbf{D}
\end{equation*}
between categories $\mathbf{C}$ and $\mathbf{D}$ is a mapping of objects to objects and arrows to arrows, in such a way that
\begin{itemize}
\item[(a)] $F(f\colon A \longrightarrow B)=F(f)\colon F(A) \longrightarrow F(B)$,
\item[(b)] $F(1_{A})=1_{F(A)}$,
\item[(c)] $F(g\circ f) = F(g)\circ F(f)$.
\end{itemize}
That is, $F$ respects domains and codomains, identity arrows, and composition.
\end{Definition}
\begin{Definition}
For categories $\mathbf{C}$, $\mathbf{D}$ and functors $F, G\colon \mathbf{C} \longrightarrow \mathbf{D}$ a {\itshape natural transformation} $\vartheta\colon F \longrightarrow G$ is a family of arrows in $\mathbf{D}$
\begin{equation*}
\left(\vartheta_{C}\colon FC \longrightarrow GC \right)_{C\in\mathbf{C}},
\end{equation*}
such that, for any $f\colon C\longrightarrow C'$ in $\mathbf{C}$, one has $\vartheta_{C'}\circ F(f) = G(f)\circ \vartheta_{C}$, that is, the following diagram commutes:
\begin{equation*}
\begin{tikzpicture}[scale=1,cap=round]
  \tikzstyle{information text}=[rounded corners,fill=red!10,inner sep=1ex]
  \begin{scope}[>=latex]
    \path (0,0) node (a) {$FC$};
    \path (2,0) node (b) {$GC$};
    \path (0,-2) node (c) {$FC'$};
    \path (2,-2) node (d) {$GC'$};
    \draw[->,line width=1.0pt] (a) -- (b) node[midway,above] {$\vartheta_{C}$};
    \draw[->,line width=1.0pt] (a) -- (c) node[midway,left] {$F f$};
    \draw[->,line width=1.0pt] (b) -- (d) node[midway,right] {$G f$};
    \draw[->,line width=1.0pt] (c) -- (d) node[midway,below] {$\vartheta_{C'}$};
  \end{scope}
\end{tikzpicture}
\end{equation*}
\end{Definition}
\begin{Definition}\label{Definition:initial_terminal_objects}
In any category $\mathbf{C}$, and object
\begin{itemize}
\item $0$ is {\itshape initial} if for any object $C$ there is a unique morphism $0\longrightarrow C$,
\item $1$ is {\itshape terminal} if for any object $C$ there is a unique morphism $C\longrightarrow 1$.
\end{itemize}
\end{Definition}
\begin{Definition}\label{Definition:subcategory}
A {\itshape subcategory} $\mathbf{S}$ of a category $\mathbf{C}$ is a collection of some of the objects and some of the arrows of $\mathbf{C}$, which includes with each arrow $f$ both the object $\mathrm{dom} f$ and the object $\mathrm{cod} f$, with each object $s$ its identity arrow $1_{\mathbf{S}}$ and with each pair of composable arrows $s\longrightarrow s' \longrightarrow s''$ their composite.
\end{Definition}
%
%\begin{specialtable}[H] 
%%\tablesize{\scriptsize}
%\caption{This is a table caption. Tables should be placed in the main text near to the first time they are~cited.\label{tab2}}
%%\tablesize{} % You can specify the fontsize here, e.g., \tablesize{\footnotesize}. If commented out \small will be used.
%\begin{tabular}{ccc}
%\toprule
%\textbf{Title 1}	& \textbf{Title 2}	& \textbf{Title 3}\\
%\midrule
%Entry 1		& Data			& Data\\
%Entry 2		& Data			& Data\\
%\bottomrule
%\end{tabular}
%\end{specialtable}
%
%\section{}
%All appendix sections must be cited in the main text. In the appendices, Figures, Tables, etc. should be labeled, starting with ``A''---e.g., Figure A1, Figure A2, etc. 

%%%%%%%%%%%%%%%%%%%%%%%%%%%%%%%%%%%%%%%%%%
\end{paracol}
%%%%%%%%%%%%%%%%%%%%%%%%%%%%%%%%%%%%%%%%%%
% To add notes in main text, please use \endnote{} and un-comment the codes below.
%\begin{adjustwidth}{-5.0cm}{0cm}
%\printendnotes[custom]
%\end{adjustwidth}
%%%%%%%%%%%%%%%%%%%%%%%%%%%%%%%%%%%%%%%%%%
\reftitle{References}

\end{document}